\documentclass{amsart}
\usepackage{amssymb,amsthm}
\usepackage{color}
\usepackage{graphicx,graphics,epsf}
\usepackage{verbatim}
\usepackage{pstricks}
\usepackage{epsfig}
\usepackage{multicol, amsmath, amssymb, amsthm,epsf,graphicx,graphics}
\usepackage{color}
\usepackage{graphicx}
\usepackage{verbatim}
\usepackage{pstricks}
\usepackage{cite}
\newtheorem{thm}{Theorem}[section]
\newtheorem{cor}[thm]{Corollary}
\newtheorem{lem}[thm]{Lemma}

\theoremstyle{definition}

\theoremstyle{remark}

\numberwithin{equation}{section}
\newcommand{\set}[1]{\left\{#1\right\}}
\newcommand{\Real}{\mathbb R}

\newcommand{\func}[1]{\ensuremath{\mathrm{#1} \:} }

\newcommand{\dist}[0]{\mathrm{dist}}
\newcommand{\re}[0]{\func{Re}}
\newcommand{\im}[0]{\func{Im}}

\newcommand{\pTorus}[1]{\mathbb{T}_{#1}^{2*}}
\newcommand{\fF}[0]{\mathbf{F}}
\title{Symmetry of Embedded Genus-one Helicoids}
\author{Jacob Bernstein and Christine Breiner}
\address{Dept. of Math,
Stanford University, Stanford, CA 94305, USA}

\email{jbern@math.stanford.edu}
\address{Dept. of Math, Massachusetts Institute of
Technology, Cambridge, MA  02139, USA} \email{breiner@math.mit.edu}
\thanks{The authors were supported respectively by the NSF grants DMS-0902721 and DMS-0902718.}

\begin{document}
\begin{abstract}
In this note, we use the Lopez-Ros deformation introduced in
\cite{LopezRos} to show that any embedded genus-one helicoid must be
symmetric with respect to rotation by $180^\circ$ around a normal
line.  This partially answers a conjecture of Bobenko from
\cite{Bobenko}.  We also show this symmetry holds for an embedded genus-$k$ helicoid $\Sigma$, provided the underlying conformal structure of $\Sigma$ is hyperelliptic.
\end{abstract}

\maketitle
 In \cite{Bobenko}, Bobenko conjectures that any immersed
genus-$k$ helicoid (i.e. a minimally immersed, once punctured
genus-$k$ surface with ``helicoid-like" behavior at the puncture) is
symmetric with respect to rotation by $180^\circ$ around a line
perpendicular to the surface. This conjecture is motivated by the observation in  \cite{Bobenko} that
the period problem for these surfaces is algebraically ``well-posed" when there is such a symmetry, but is ``over-determined'' without it. In this note, we verify Bobenko's conjecture for
\emph{embedded} genus-one helicoids.  That is:

\begin{thm} \label{RotSymThm}
Let $\Sigma$ be an embedded genus-one helicoid.  Then there is a line $\ell$ normal to $\Sigma$
so that
% $\Sigma$ is symmetric with respect to
%rotation by $180^\circ$ around $\ell$.
rotation by $180^\circ$ about $\ell$ acts as an orientation preserving isometry on $\Sigma$.

\end{thm}

We define a genus-$k$ helicoid to be a complete, minimal surface
immersed in $\Real^3$ which has genus $k$, one end, and is
asymptotic to a helicoid.   A consequence of Theorem 3 of \cite{HPR}
is that any
(immersed) minimal surface which is conformally a once-punctured
compact genus-$k$ Riemann surface with ``helicoid-like" Weierstrass
data at the puncture is a genus-$k$ helicoid in this sense.  In particular, the above
definition encompasses the surfaces studied by Bobenko. Importantly,
by Theorem 1.1 of \cite{BB2}, any complete, \emph{embedded} minimal
surface in $\Real^3$ with genus $k$ and one end has ``helicoid-like" Weierstrass data and hence is a genus-$k$
helicoid.  The space of such objects is not vacuous.  Weber, Hoffman and Wolf \cite{WHW} and Hoffman and White
\cite{HW} have given (very different) constructions of embedded
genus-one helicoids --  at present it is unknown whether the two constructions give the same surface.  Both constructions produce a genus-one helicoid that has, in addition to the orientation preserving symmetry of Theorem \ref{RotSymThm}, two orientation reversing symmetries.  Whether all genus-one helicoids possess these additional symmetries is also unknown.

We emphasize that our argument does not generalize to genus $k>1$
because we crucially use the fact that every genus-one Riemann
surface admits a large number of biholomorphic involutions -- more
precisely, that any once-punctured genus-one Riemann surface admits a
non-trivial biholomorphic involution. This need not be true for higher
genus. Indeed, \emph{a priori} there may fail to be any non-trivial
biholomorphic automorphisms.  However, if we restrict attention to
genus-$k$ helicoids whose underlying Riemann surface structure is
\emph{hyperelliptic} -- that is the surface admits a biholomorphic
involution $I$ with $2k+2$ fixed points, one of which is the
puncture, then our arguments continue to hold.  Consequently, we present the argument in this more general context.  Finally, we note that Francisco Martin has pointed out to us that 
with only slight modifications, the argument also proves that embedded \emph{periodic} genus-one helicoids admit such a symmetry.  

Let us outline the proof of Theorem \ref{RotSymThm} for genus-one helicoids.  By Theorem 1.1
of \cite{BB2}, $\Sigma$ is conformally a once-punctured torus with
``helicoid-like" Weierstrass data at the puncture.  Thus, $\Sigma$
admits a biholomorphic involution, $I$, which is compatible with
this data.  Indeed, if $dh$ is the height differential and $g$ is
the stereographic projection of the Gauss map then $I^*dh=-dh$ and
$g\circ I=C g^{-1}$, for $C\in \mathbb{C}\backslash \set{0}$.  If
$|C|= 1$, a simple computation using the Weierstrass representation
implies Theorem \ref{RotSymThm}. On the other hand, if $|C|\neq  1$
then the interaction between the period conditions and the
involution $I$ imply $\Sigma$ has vertical flux. In this case,
following Perez and Ros \cite{PerezRos}, we may deform the
Weierstrass data to obtain  a smooth family of immersed minimal
surfaces, $\Sigma_\lambda$.  Here  $\Sigma=\Sigma_1$ and $\Sigma_\lambda$ is the Lopez-Ros deformation \cite{LopezRos} of $\Sigma$.  As in
\cite{PerezRos}, for $\lambda$ near $1$, $\Sigma_\lambda$ is
embedded, while for $\lambda>>1 $, $\Sigma_\lambda$ is not embedded,
contradicting the maximum principle.

%Indeed, by Theorem 3 of \cite{HPR}, each $\Sigma_\lambda$ is a genus-one helicoid and hence has an embedded end. Thus,  for $\lambda$ close to $1$, $\Sigma_\lambda$ is embedded. On the other hand, as $\Sigma$ must contain a point with vertical normal, Lemma 4 of \cite{PerezRos} implies that for $\lambda$ large, $\Sigma_\lambda$ has self-intersection (i.e. is not embedded). Combining these, there must be a first $\lambda>1$ for which $\Sigma_\lambda$ has a  self-intersection and any self-intersection point lies within  a compact set.  This contradicts the strong maximum principle. The argument we use is inspired by Perez and Ros's proof in \cite{PerezRos} of the result of Lopez and Ros \cite{LopezRos} on the non-existence of complete, embedded, minimal planar domains with finite total curvature and more than 2 ends.

\subsection*{Acknowledgments}
We would like to thank David Hoffman, Francisco Martin and the
anonymous referees for their many constructive comments.  We
also thank Brian White for clarifying some
results regarding mass minimizing currents.
\section{Asymptotic properties of $\Sigma$ and properties of the Involution}
\subsection{The Weierstrass Representation and the Flux} \label{WeierstrassSec}
We recall the Weierstrass representation for immersed
minimal surfaces in $\Real^3$.  Let $M$ be a Riemann surface and
suppose that $g$ is a meromorphic function on $M$ and $dh$ a
holomorophic one-form.  Suppose, moreover, that the meromorphic one-forms $g dh$ and $g^{-1} dh$ have no poles and do not simultaneously vanish.  Then the map $\fF:
{M} \to \Real^3$ given by
\begin{equation}\label{WeierstrassRep}
(x_1,x_2,x_3)=\mathbf{F}:=\re \int\left(\frac{1}{2}
(g^{-1}-g),\frac{i}{2}(g^{-1}+g), 1\right) dh
\end{equation}
is a minimal immersion with the property that $g$ is the
stereographic projection of the Gauss map of the image of
$\mathbf{F}$ and $\re dh=\mathbf{F}^*dx_3$. Without further
restrictions on the data, $\fF$ is potentially only defined on
$\tilde{M}$, the universal cover of $M$. These restrictions are
known as the \emph{period conditions}, which, when satisfied, ensure
that $\fF$ is well-defined on $M$.  Explicitly they may be stated
as:
\begin{equation} \label{PeriodCond}
 \int_{\gamma} g dh=\overline{\int_\gamma g^{-1}dh } \; \; \mbox{and} \; \; \re \int_{\gamma} dh=0
\end{equation}
for any closed curve, $\gamma$, on $M$.  Conversely, given a minimal immersion
$\fF:M\to \Real^3$, one obtains $g$ and $dh$ satisfying the period
conditions and with $g dh, g^{-1} dh$ holomorphic and not vanishing simultaneously and so that the image of the map given by
\eqref{WeierstrassRep} coincides with the image of $\fF$ (up to a
translation).

We will also consider the \emph{flux} of the immersion $\fF$ along
closed curves.  For $\gamma$ a closed curve on $\Sigma$, we denote by
$\mathbf{\nu} = -d\fF (J \gamma')$ the conormal vector field along $\gamma$.
Here $J$ is the complex structure of $\Sigma$ and $\gamma'$ is the
derivative of $\gamma$ with respect to arc length.  We define the flux
of $\Sigma$ along $\gamma$ equivalently as:
\begin{equation}\label{FluxDef}
Flux(\gamma)= \int_\gamma \mathbf{\nu}\; ds= \im \int_\gamma \left(\frac{1}{2}
(g^{-1}-g),\frac{i}{2}(g^{-1}+g), 1\right) dh.
\end{equation}
The equivalence of the two definitions is a simply consequence of the Cauchy-Riemann equations and \eqref{WeierstrassRep}.  Indeed, on an oriented Riemannian surface, every harmonic one-form $\omega$ has a harmonic conjugate $\omega^*=-\omega\circ J$ with $\omega+i \omega ^*$ holomorphic.  Conversely,  any holomorphic one-form can be written as $\omega+i \omega^*$ with $\omega$ harmonic.
As holomorphic one-forms are closed,  Stokes' theorem implies the
flux of a curve depends only on its homology class.

Of particular interest in this paper are surfaces with everywhere
\emph{vertical flux}, i.e. where the horizontal components of the
flux are zero for all closed curves. A
minimal surface with vertical flux can be smoothly deformed to give
a smooth family of minimal immersions.
Indeed, suppose $\Sigma$ is a minimal surface with vertical flux and
Weierstrass data $(g,dh,M )$. Then for $\lambda \in \Real^+$, it
follows from \eqref{FluxDef} that the
triple $(\lambda g, dh, M)$ satisfies the period conditions
\eqref{PeriodCond} and so \eqref{WeierstrassRep} gives the desired family of immersions
$\fF_\lambda:M\to \Real^3$. 
Such a deformation was introduced by Lopez and Ros in \cite{LopezRos}.
%We now focus on the case $\Sigma \in
%\mathcal{E}(1,1)$.  By Theorem 1.1 with genus one as $\Sigma$ admits a
%holomorphic involution, it acutally has a rotational symmetry. Indeed by studying the action of this involution
%on the meromorphic one-form data we conclude that $\Sigma$  a rotation around the $x_1$-axis, thus proving
%Corollary \ref{AppCor}.
\subsection{The Involution of $\Sigma$}
We now consider $\Sigma$, an embedded genus-$k$ helicoid with
asymptotic helicoid $H$. Denote the Weierstrass data of $\Sigma$ by
$(g,dh,M)$.  By Theorem 1.1 of \cite{BB2},  $\Sigma$ is conformal to
a once-punctured compact genus-$k$ surface and the one-forms $dh$ and $\frac{dg}{g}$
both have double poles at this puncture with zero residue -- this is the ``helicoid-like" behavior alluded to in the introduction.  
More precisely,  there is a compact Riemann surface $M_k$ and a point
$\infty\in M_k$ so that $M=M_k\backslash \set{\infty}$ with $dh$ and
$\frac{dg}{g}$ meromorphic one forms on $M_k$ both with a double pole
at $\infty$ and no residue there.  We assume also that $\Sigma$ is
\emph{hyperelliptic}.  That is, there exists a non-trivial
biholomorphic involution $I:M_k\to M_k$ with $2k+2$ fixed points and
so that $I(\infty)=\infty$.  An important property of hyperelliptic
involutions is that
for $0\neq [\gamma]\in H_1(M_k, \mathbb{Z})$ a non-trivial element of the first
homology group of $M_k$, $I_*[\gamma]=-[\gamma]$ -- see
\cite{FarkasKra}. As inclusion of $M$ in $M_k$ induces an isomorphism
between $H_1(M,\mathbb{Z})$ and $H_1(M_k, \mathbb{Z})$ this property also holds for $\Sigma$.
 
Hyperellipticity is a very strong condition when
$k>1$.  
 However, genus-one helicoids are always hyperelliptic. Indeed, let $\Lambda_\tau=\lbrace n+m\tau : n,m\in
\mathbb{Z}\rbrace \subset\mathbb{C}$  be the lattice so that
$\pTorus{\tau}=\mathbb{C}/\Lambda_\tau\backslash \set{\bar{0}}$ is conformally
equivalent to $\Sigma$, where $\bar{0}=0+\Lambda_\tau$.  As
$-\Lambda_\tau=\Lambda_\tau$, the map $u \to -u$ induces a
biholomorphic involution of $\pTorus{\tau}$
%$\mathbb{C}/\Lambda_\tau\backslash \bar{0}$
and
hence an involution $I:\Sigma \to \Sigma$. The
half-period lattice of $\Lambda_\tau$ is fixed by $u\to -u$, and so $\infty$
and exactly 3 points of $\Sigma$ are fixed by $I$.

Before we proceed we note the following simple lemma:
\begin{lem} \label{CpxVarLem} 
Suppose that $p$ is a point in $N$, a Riemann surface with a non-trival involution $I:N\to N$, so that $I(p)=p$. Then there is a coordinate neighborhood $U$ about $p$ with coordinate $u$ so that: $I(U)=U$
and $u\circ I=-u$. 
Moreover, suppose $\omega$ is a meromorphic one-form on $U$ of the form:
\begin{equation*}
\omega=\left(\frac{a}{u^2}+\frac{b}{u}+H(u)\right) du
\end{equation*}
with $H$ holomorphic.  Then the one-form $\omega+I^*\omega$
\begin{enumerate}
\item Has a simple pole at $p$ iff $b\neq 0$,
\item Has a zero at $p$ iff $b=0$,
\item \label{vanishing}Is identically zero on $U$ iff $b=0$ and $H$ is even.
\end{enumerate}
\end{lem}
By \eqref{vanishing}, if $H$ has a simple zero at $u=0$ then $\omega +I^* \omega$ can't vanish identically on $U$.
\begin{proof}
The existence of such a coordinate is a straightforward consequence of the inverse function theorem.  One calculates that in $U$,  $I^* \omega=-\left(\frac{a}{u^2}-\frac{b}{u}+H(-u)\right) du$ and so 
$\omega+I^* \omega=\left( \frac{2b}{u} +H(u)-H(-u)\right) du$. Clearly
$H(u)-H(-u)$ is odd and so vanishes at $u=0$; from this all three claims follow. 
\end{proof}

We now analyze how $I$ acts on the Weierstrass data:
\begin{lem}\label{InvOndh} $I^* dh=-dh$ and $I^* \frac{dg}{g}=-\frac{dg}{g}$
\end{lem}
\begin{proof}
Let us denote by $F\subset M_k$ the fixed point set of $I$. For a given meromorphic one-form $\alpha$, let $Z[\alpha]$ and $P[\alpha]$ represent the sets, respectively,  of zeros and poles of $\alpha$.
The Riemann-Roch theorem implies that for any non-vanishing meromorphic one-form $\alpha$ on $M_k$ with $Z[\alpha]\neq \emptyset$  one has the following relation:
\begin{equation}\label{RiemannRochEqn}
\# Z[\alpha]-\# P[\alpha]=2k-2
\end{equation}
where $\# Z[\alpha]$ and $\# P[\alpha]$ denote the number of, respectively, zeros and poles of $\alpha$ counting multiplicity.  For an arbitrary set of points, $X$, we denote by $|X|$ the number of points of $X$.  In general, $\# Z[\alpha]\geq |Z[\alpha]|$ and $\# P[\alpha]\geq |P[\alpha]|$. 

Recall that $dh$ has only one pole (at $\infty$) and it is a double pole with no residue.  Given that $\infty$ is fixed by $I$, Lemma \ref{CpxVarLem} implies that $I^*dh+dh$ has no poles and has zeros at each point of $F$.  As $|F|=2k+2$, \eqref{RiemannRochEqn} implies that $I^*dh+dh$ must vanish identically.  This proves the first part of the lemma and that $I$ preserves $Z[dh]$.  Lemma \ref{CpxVarLem}, in particular \eqref{vanishing}, and \eqref{RiemannRochEqn} then imply that $|Z[dh]|\leq 2k-|F\cap Z[dh]|$. 

Set $\omega=\frac{dg}{g}$ and $\tilde{\omega}=\omega+I^*\omega$. In $\Sigma$, the poles of $\omega$ are simple and,  as $g dh$ and $g^{-1} dh$ do not simultaneously vanish,  these poles occur precisely at the points of $Z[dh]$.   Thus, Lemma \ref{CpxVarLem} implies that: $\tilde{\omega}$ has only simple poles; $P[\tilde{\omega}]\subset Z[dh]$; and $Z[\tilde{\omega}]\supset F\backslash (F\cap Z[dh])$.  
As the poles are simple:
\begin{equation}\label{poletilde}
\# P[\tilde{\omega}]\leq |Z[dh]|\leq  2k-|F\cap Z[dh]|
 \end{equation}
while
\begin{equation}\label{zerotilde}
\# Z[\tilde{\omega}]\geq |F|-|F\cap Z[dh]|=2k+2-|F\cap Z[dh]|.
\end{equation}
If $\tilde{\omega}$ does not vanish identically then \eqref{RiemannRochEqn} implies $k\geq 2$.  This proves the theorem for $k=1$.   

For $k>1$ we must use further properties of genus-$k$ helicoids.  To begin the argument assume $\tilde \omega$ is not identically zero.  Observe that in $\Sigma$, $\omega$ has only simple poles which occur at the zeros and poles of $g$.  The residue of $\omega$ at such a zero or pole is exactly equal to $\pm m$ where $m$ is the order of the zero or pole.  
%Note that the magnitude of the residue corresponds to the order of vanishing of $dh$.  
 Lemma \ref{AppendixLem} proves that $p$ is a pole of $g$ if and only if $I(p)$ is a zero of the same order.  
 %Since we know $I$ preserves the order of vanishing of $dh$, 
 Thus, the residues of ${\omega}$ and of $I^* \omega$ cancel at any
 pole of $\omega$ in $\Sigma$.  Hence, $\#P[\tilde{\omega}]=0$ and by
 \eqref{RiemannRochEqn} if $\tilde{\omega}$ doesn't vanish identically
 then $\#Z[\tilde \omega]=2k-2$.  
 
 Finally, Lemma \ref{AppendixLem} implies that if $p$ is a zero or pole of $g$, then $p \notin F$.  As $P[\omega] = Z[dh]$, it follows that $F \cap Z[dh]= \emptyset$.  Then by \eqref{zerotilde}, $\#Z[\tilde \omega] \geq 2k+2$.  This gives the necessary contradiction and completes the proof.
\end{proof}
We next compute how the Gauss map is transformed under
$I$. First, pick $p_0 \in\Sigma$ satisfying
$I(p_0)=p_0$ so that $\frac{dg}{g}$ does not have a pole at $p_0$
-- there are $2k+1$ points of $\Sigma$ fixed by $I$ and by \eqref{RiemannRochEqn} at most $2k$ poles of
$\frac{dg}{g}$, so such a point $p_0$ exists. We determine the transformation of the
Gauss map using its value at this fixed point:
\begin{cor}\label{InvOng}
$g\circ I=\frac{g(p_0)^2}{g}$ for $g(p_0)\in
\mathbb{C}\backslash\set{0}$ and $p_0$ as determined above.
\end{cor}
\begin{proof}
By analytic continuation it suffices to consider $U$ the neighborhood of $p_0$ from Lemma \ref{CpxVarLem}.  We have:
$\frac{g(p)}{g(p_0)}=\exp(\int_\gamma \frac{dg}{g})$ where $\gamma$ is a
path in $U$ connecting $p_0$ to $p$.  Then $\frac{g(I(p))}{g(p_0)}=\exp(\int_{I(\gamma)}
\frac{dg}{g})$. However, $\int_{I(\gamma)} \frac{dg}{g}=\int_\gamma
I^*\frac{dg}{g}=-\int_\gamma \frac{dg}{g}$ and so
$g(I(p))=\frac{g(p_0)^2}{g(p)}$.  
\end{proof}

\section{The Rotational Symmetry} \label{RotSymSec}
Using the properties of the involution $I$, we can now prove that
$\Sigma$ has the claimed symmetry.  Note that by rotating $\Sigma$
about the $x_3$-axis and translating $\Real^3$, we may assume that
$g(p_0)>0$ and $\fF(p_0)=0$; here $p_0$ is the point from Corollary
\ref{InvOng}.  If $g(p_0)=1$ then a simple computation using the
Weierstrass representation gives that $(x_1,x_2,x_3)\circ I=(x_1,
-x_2, -x_3)$,  proving Theorem \ref{RotSymThm}. Thus, we
must rule out the possibility that $g(p_0)\neq 1$.

To that end, we use $I$ to see that in this
case $\Sigma$ has vertical flux:
\begin{lem}\label{gneq1lem}
If $g(p_0)\neq 1$ then $g dh$ and $\frac{1}{g} dh$ are exact forms
on $\Sigma$.
\end{lem}
\begin{proof}
Recall $g dh$ and $\frac{1}{g} dh$ are both holomorphic one-forms on
$\Sigma$ and are hence closed.  As a consequence, it will suffice to
show that over the $2k$ generators of the homology group $[\eta_i]$ that $\int_{\eta_i} g dh=\int_{\eta_i} \frac{1}{g} dh =
0$. Here $\eta_i$ are simple closed curves and $[\eta_i]$ the
corresponding homology classes. For simplicity, we treat only $g
dh$. By the first equation in
\eqref{PeriodCond}:
\begin{equation*}
\int_{\eta_i} g dh =\overline{\int_{\eta_i} \frac{1}{g} dh}.
\end{equation*}
Recall, hyperelliptic involutions satisfy $I_*[\eta_i]=-[\eta_i]$. Hence, $\int_{\eta_i} g dh
=-\int_{I(\eta_i)} g dh=g(p_0)^2 \int_{\eta_i} \frac{1}{g}
dh=g(p_0)^2 \overline{\int_{\eta_i} g dh}$. Taking absolute values,
if $g(p_0) \neq 1$, then $\int_{\eta_i} g dh
=0$.
\end{proof}

We will argue as in \cite{PerezRos} to show that the existence of a
$\Sigma$ with vertical flux is precluded by the maximum principle.
First we show:
\begin{lem} \label{VertFluxConLem}
Suppose $\Sigma$ has vertical flux.  Then, there is a smooth family of immersed minimal surfaces $\Sigma_\lambda$, $\lambda>0$, with $\Sigma_1=\Sigma$ and a fixed helicoid $H$ so that:
\begin{enumerate}
\item \label{itemi} Each $\Sigma_\lambda$ is a genus-$k$ helicoid and is asymptotic to $H$.
\item\label{itemii} The set $E=\set{\lambda\in \Real^+ : \Sigma_\lambda \mbox{ is embedded}}$ is open.
\item \label{itemiii} For $\lambda$ sufficiently large $\Sigma_\lambda$ is not embedded.
\end{enumerate}
\end{lem}
\begin{proof}
As $\Sigma$ has vertical flux, the triple $(\lambda g, dh, M)$, for
$\lambda \in  \Real^+$,  gives rise to a minimal immersion
$\fF_\lambda:M_k\to \Real^3$.  Let us denote by $\Sigma_\lambda$ the
image of $\fF_\lambda$ and  set $g_\lambda= \lambda g$. Notice that
each $\Sigma_\lambda$ is a complete, minimally immersed genus-$k$ surface and
$\Sigma_1=\Sigma$.  It remains only to verify that this family
satisfies \eqref{itemi}-\eqref{itemiii}.

To that end, we note that by Corollary 1.2 of \cite{BB2}, $M_k$
has a neighborhood of infinity, $U$, with holomorphic coordinate
$z:U\backslash \set{\infty}\to \mathbb{C}$ so that (after possibly rescaling $\Sigma$) on
$U\backslash \set{\infty}$, $g(p)=\exp( iz(p)+F(p))$, where $F$ is holomorphic and has a
zero at $\infty$.  As a consequence,
$g_\lambda(p)=\exp(iz(p)+\log \lambda+ F(p))$.  Hence, Theorem 3 of
\cite{HPR} implies  that,
outside of a ball of radius $R_\lambda$, $\Sigma_\lambda$ is
asymptotic to a scale 1 helicoid -- i.e. a helicoid with Weierstrass
data $(e^{iz},dz, \mathbb{C})$. After a translation, this gives \eqref{itemi}. To see
\eqref{itemii}, we note that for any $\lambda_0$ there is an $R>1$
so that, for $\lambda\in [\lambda_0/2,2\lambda_0]$, outside of $B_R$
each $\Sigma_\lambda$ is a normal exponential graph over the
helicoid $H$ with small $L^\infty$ norm.  In particular, the
$\Sigma_\lambda$ are embedded outside $B_R$.  On the other hand as
$\lambda \to \lambda_0$ the $B_R\cap \Sigma_\lambda$ converge
smoothly to $B_R\cap \Sigma_{\lambda_0}$. By \eqref{itemi} this
convergence must be with multiplicity 1 and so for $\lambda$ close
enough to $\lambda_0$, $B_R\cap \Sigma_\lambda$ can be written as
the normal exponential graph with small $L^\infty$ norm over
$B_R\cap \Sigma_{\lambda_0}$. Thus, if $B_{R}\cap
\Sigma_{\lambda_0}$ is embedded then, for $\lambda$ sufficiently
close to $\lambda_0$, so is $B_R\cap \Sigma_\lambda$. This gives
\eqref{itemii}. Finally, \eqref{RiemannRochEqn} implies that $g$ must have at least one pole and one zero
on $\Sigma$. Lemma 4 of \cite{PerezRos} then gives \eqref{itemiii}.
\end{proof}

\begin{proof} (of Theorem \ref{RotSymThm})
By the above, it suffices  to show that $\Sigma$ does not have vertical flux.
We proceed by contradiction. If $\Sigma$ has vertical flux then
Lemma \ref{VertFluxConLem} gives a family $\Sigma_\lambda$ with
the properties \eqref{itemi}-\eqref{itemiii}.  By \eqref{itemii} and
\eqref{itemiii}, there exists a $\lambda_0>1$ so that
$\Sigma_{\lambda_0}$ is not embedded, but for $\lambda\in
[1,\lambda_0)$, $\Sigma_\lambda$ is embedded.  By \eqref{itemi},
\eqref{itemiii}, and the fact that the $\Sigma_\lambda$ smoothly
depend on $\lambda$, there are points $p_\lambda^1, p_\lambda^2$
with $|p_\lambda^i|\leq R<\infty$ so that $\lim_{\lambda \nearrow
\lambda_0} |p_\lambda^1-p_\lambda^2|\to 0$ but
$\dist_{\Sigma_\lambda} (p_\lambda^1,p_\lambda^2)\geq \delta>0$. By
the strong maximum principle, this is only possible if the
$\Sigma_\lambda$ converge to $\Sigma_{\lambda_0}$ as $\lambda\nearrow
\lambda_0$ with multiplicity greater than 1. This contradicts
\eqref{itemi}. 
\end{proof}
\section{Periodic Genus-One Helicoids}
Francisco Martin has kindly pointed out to us that our argument can be
readily adapted to embedded \emph{periodic} genus-one helicoids.  For
the sake of completeness, we
include here a sketch of his argument.   

Roughly speaking, a periodic genus-one helicoid looks like  a helicoid  with an infinite number of handles placed periodically along the axis.  More precisely, it
 is an infinite genus surface, $\Sigma$ so that: $\Sigma$ is asymptotic to a helicoid; $\Sigma$ is invariant under a ``screw-motion" $S_{\theta,t}$; and the quotient of $\Sigma$ by the group generated by $S_{\theta,t}$ is a twice punctured torus.  Here $S_{\theta,t}$ is the isometry of $\Real^3$ given by rotating about the $x_3$-axis by $\theta$ followed by a translation by $t$ in the $x_3$ direction.  
 There exists a family of embedded examples.  Indeed, in \cite{WHW} an embedded genus-one helicoid is constructed as a limit of periodic genus-one helicoids. See also \cite{HKW}.  

As discussed in \cite{AFM}, after a homothety, one may consider a periodic genus-one helicoid to be a triple of Weierstrass data:
$(g, dh,\mathbb{T}^2 \backslash \set{E_1, E_2})$.  Here $E_1,E_2$ are distinct points of $\mathbb{T}^2$ and $g$ extends meromorphically to $\mathbb{T}^2$ with a zero at $E_1$ and a pole at $E_2$.  Further, $dh$ extends meromorphically to $\mathbb{T}^2$ with simple poles at $E_1, E_2$, and residue $-i$ at $E_1$ and $i$ at $E_2$. Finally, the compatibility conditions of Section \ref{WeierstrassSec} hold except that the vertical periods around $E_1$ and $E_2$ do not close up; recall we are parameterizing a surface of infinite topology.  For simplicity, take this as our definition and refer to \cite{AFM, MRPer} for weaker characterizations.
 It is a straightforward computation to see that a periodic genus-one helicoid is asymptotic to a scale-one helicoid if and only if  the one-form $\frac{dg}{g}-i dh$ has no poles at $E_1,E_2$ -- see \cite{MRPer, PerezRos}.  
 
\begin{lem} \label{PerLem}
There is a non-trivial biholomorphic involution $I$ of $\mathbb{T}^2$ so that $I(E_1)=E_2$.  Moreover, $I^*dh=-dh$ and $I^* \frac{dg}{g}=-\frac{dg}{g}$.
\end{lem}
\begin{proof}
Identify $\mathbb{T}^2$ with the quotient $\mathbb{C}/\Lambda_\tau$ where $\Lambda_\tau=\set{n+m \tau: n,m\in \mathbb{Z}}$.  As translation along $1$ or $\tau$ in $\mathbb{C}$ induce  biholomorphic automorphisms, we may represent $E_i$ by points $p_i+\Lambda_\tau$ where the $p_i$ are placed symmetrically with respect to $0$. Hence, the map $u\to -u$ on $\mathbb{C}$ descends to an involution $I$ of $\mathbb{T}^2$ that swaps $E_1$ and $E_2$.  As $I$ swaps the $E_i$ and the residues of $dh$ at $E_1$ and $E_2$ are of opposite sign, $I^*dh+dh$  has no poles.  By Lemma \ref{CpxVarLem}, $I^*dh+dh$ has at least four zeros and so by \eqref{RiemannRochEqn} vanishes identically.  

By construction, $\frac{dg}{g}$ has simple poles at $E_1$ and $E_2$ with residues of opposite sign.
Thus, $I^*\frac{dg}{g}+\frac{dg}{g}$ has no residue at either $E_1$ or $E_2$ and hence no poles there.  On the other hand, all other poles of $\frac{dg}{g}$ occur at the  zeros of $dh$ and these are involuted by $I$ and so $I^*\frac{dg}{g}+\frac{dg}{g}$ has at most two poles. By Lemma \ref{CpxVarLem}, this form has at least three zeros and so by \eqref{RiemannRochEqn} must vanish identically.  
\end{proof}
\begin{cor}
Let $\Sigma $ be an embedded periodic genus-one helicoid.
 Then there is a line $\ell$ normal to $\Sigma$
so that
% $\Sigma$ is symmetric with respect to
%rotation by $180^\circ$ around $\ell$.
rotation by $180^\circ$ about $\ell$ acts as an orientation preserving isometry on $\Sigma$.
\end{cor}
\begin{proof}
The corollary follows from Lemma \ref{PerLem} and the arguments of Section \ref{RotSymSec} as long as we can rule out the existence of an embedded periodic genus-one helicoid with vertical flux.  Notice that, by construction, the periods around $E_1$ and $E_2$ always have vertical flux. 
 Suppose $\Sigma$ is a periodic genus-one helicoid with vertical flux,
 asymptotic to some helicoid $H$. Let $\Sigma_\lambda$ denote the
 family of surfaces given by the Lopez-Ros deformation.  Necessarily,
 this family remains in the class of periodic genus-one helicoids and
 all have the same asymptotic behavior.  Thus, outside of a bounded
 cylinder, each $\Sigma_\lambda$ is embedded and asymptotic to $H$.
 Due to the periodicity, the non-compactness of the
 cylinder does not introduce additional difficulties.
 Clearly, \eqref{RiemannRochEqn} implies $g$ has a pole or zero in
 $\mathbb{T}^2\backslash \set{E_1,E_2}$ and so for $\lambda>>1$,
 $\Sigma_\lambda$ fails to be embedded.  Hence, the $\Sigma_\lambda$
 satisfy the conclusions of Lemma \ref{VertFluxConLem} and so one
 obtains a contradiction exactly as in the proof of Theorem
 \ref{RotSymThm}.  
 %Moreover, as $\Sigma$ is embedded $\frac{dg}{g}-i dh$ has no poles at $E_1,E_2$ and this is preserved for all $\Sigma_\lambda$. So the end of each $\Sigma_\lambda$ is asymptotic to a helicoid and hence embedded.
% The embeddedness of each $\Sigma_\lambda$, outside of a compact set, follows from the embeddedness of $\Sigma$.  Indeed, by \cite{MRPer}, we write the end of $\Sigma$ near $E_1$ in conformal coordinates on a punctured disk $D^*$.  Then $x_3(z) = arg(z) + f_1(z)$ and $g(z)= zf_2(z)$ where $f_1$ is bounded and $f_2(0) \neq 0$ (and both $f_i$ are holomorphic on $D$).  The asymptotic behavior of level curves of $x_3$ are unaffected by the Lopez-Ros deformation, which shows that each $\Sigma_\lambda$ is, at least outside of a compact set, embedded.  Furthermore they are all asymptotic to the same helicoid $H$.  
 
 %To see this, note the asymptotic behavior of $\Sigma$ shows that $\{x_3=c\} \cap \Sigma \backslash B_R$, for some large, fixed $R$, is the union of two proper, smooth curves.  We parameterize these curves as $\gamma(\pm t)$.  
 
 \end{proof}
\appendix
\section{Hyperelliptic case}
In this appendix we complete the proof of Lemma \ref{InvOndh}.  The arguments are of a rather different flavor than the rest of the paper and are a refinement of those used in \cite{BB2} to show that $\frac{dg}{g}$ had no residue at $\infty$.  We first recall the following elementary facts about level sets of harmonic functions:
\begin{lem} \label{HarmonLevelLem}
Let $V$ be an open set in a Riemannian surface $\Sigma$ with $f$ a harmonic function on $V$.  If $p\in V$ is a critical point of $f$ and $t_0=f(p)$ then:
\begin{enumerate}
\item \label{HLLi} There is a simply-connected neighborhood $U(p)$ of $p$ so that $\set{f=t_0}\cap U(p)$ consists of $m+1$ smooth embedded curves $\sigma_i$ with $\partial \sigma_i \subset \partial U(p)$;  $m$ is the order of vanishing of $f$ at $p$. The $\sigma_i$ meet only at $p$ and do so transversally. 
\item \label{HLLii}There is a decomposition $\set{f>t_0}\cap U(p)=\Sigma^+(p)=\Sigma^+_1(p) \cup \ldots \cup \Sigma^+_{m+1}(p)$ and $\set{f<t_0}\cap U(p)=\Sigma^-(p)=\Sigma^-_1(p) \cup \ldots \cup \Sigma^-_{m+1}(p)$ so  for $i\neq j$ the closure of ${\Sigma}_i^\pm $ meets the closure of ${\Sigma}_j^\pm $ only at $p$.   
\item \label{HLLiii}For $t>t_0$, the set $\Sigma^{\geq t}(p)=\set{f\geq  t}\cap U(p)$ consists of $m+1$ components each in a different component of $\Sigma^+(p)$. 
\item \label{HLLiv}For each $i$, there is a piecewise smooth parameterization $\gamma_i:(-1,1)\to M$ of $\partial \Sigma^\pm_i(p)\cap U(p)$ and a sequence of smooth injective maps $\gamma_i^j:(-1,1)\to \Sigma^\pm_i(p)$ so that: on $(-1,-1/2)\cup (1/2, 1)$, $\gamma_i^j=\gamma_i$ and $\gamma_i^j\to \gamma_i$ in $C^0$ as $j\to \infty$.
\end{enumerate}
\end{lem}

We also note the following facts regarding the asymptotic properties of level curves of the height function of a genus-$k$ helicoid:
\begin{lem} \label{CylLem} Given $\Sigma$ an embedded genus-$k$ helicoid, there is a cylinder:
  \begin{equation*}
  C= C_{h,R}=\set{|x_3|\leq h, x_1^2+x_2^2\leq R^2}
  \end{equation*}
and a component $\Sigma'$ of $C\cap \Sigma$ so that: 
\begin{enumerate}
\item all critical points of $x_3:\Sigma\to \Real$ lie on the interior of $\Sigma'$.
\item $\partial \Sigma'=\gamma_t\cup \gamma_b\cup \gamma_u \cup \gamma_d$, four smooth curves, with $x_3=h$ on $\gamma_t$, $x_3=-h$ on $\gamma_b$ and for $s\in (-h,h)$, $\set{x_3=s}$ meets $\gamma_u$ and $\gamma_d$ each in one point.
\end{enumerate}
\end{lem}
\begin{proof}
As $\Sigma$ is properly embedded there
exist $h$ and $R$ so all the zeros and poles of $g$ lie in the  interior of the 
cylinder $C$.
Moreover, by increasing the size of the cylinder one can take a component,  $\Sigma'$, of $\Sigma \cap C$ so that all of the zeros and poles of $g$ lie in $\Sigma'$ and, as $\Sigma$ has one end, so that $\Sigma\backslash \Sigma'$ is an annulus. 
Finally, as $\Sigma $ is asymptotic to some helicoid $H$, by further enlarging the cylinder, we may take $\Sigma'$ so $\gamma=\partial
\Sigma'$ is the union of four smooth curves, two at the top and
bottom, $\gamma_t$ and $\gamma_b$, and two disjoint helix like
curves $\gamma_u,\gamma_d$  so
$\frac{d}{dt}x_3(\gamma_u(t))>0$ and $\frac{d}{dt}x_3(\gamma_d(t))<0$.
\end{proof}

\begin{lem} \label{AppendixLem}
A point $p\in \Sigma$ is a pole of $g$ if and only if $I(p)$ is a zero of $g$ of the same order.
\end{lem}
\begin{proof}
First note that  Lemma \ref{InvOndh} implies $I^* dh=-dh$; thus, if
$p$ is a zero of order $m$ of $dh$ so is $I(p)$ and, up to
a vertical translation, $x_3\circ I=-x_3$ .  Recall, $g dh$ and $g^{-1} dh$ are holomorphic in $\Sigma$ and do not
simultaneously vanish, hence the order of a pole of $g$ or zero of $g$
 at $p$ is equal to the order of the
zero of $dh$ at $p$.  Thus, it suffices to show that if $p$ is a zero of $g$
then $I(p)$ is a pole.

Let $R, h, C$ and $\Sigma'$ be as in Lemma \ref{CylLem}.  It is a standard topological fact that a closed, oriented and connected surface in $\Real^3$ divides $\Real^3$ into two components. Thus, $C\backslash \Sigma'$ consists of two components $\Omega^+$ and $\Omega^-$, labeled so that the normal, $\mathbf{n}$, to $\Sigma$ points into $\Omega^+$.  
Denote by $\sigma_t$ the set $\Sigma'\cap \set{x_3=t}$ and by
$\Omega_t^\pm$ the set $\Omega^\pm \cap \set{x_3=t}$.  
A fact we will
use below is that the closed sets $\bar{\Omega}^\pm$ are the complements
of the union of open sets with smooth boundaries.  Indeed, let $U_1=U_1^\pm$
be the component of $\Real^3\backslash \Sigma$ containing $\Omega^\mp$, $U_2=\set{x_1^2+x_2^2>R^2}$, $U_3=\set{x_3>h}$ and
  $U_4=\set{x_3<-h}$.  Then all the $U_i$ are open with smooth boundary
  and $\bar{\Omega}^\pm=\Real^3\backslash \cup_{i=1}^4 U_i$.  

Consider a critical point $p$ of $x_3$ with $t_0=x_3(p)$ and $x_3$ vanishing to order $m$ at $p$.  At $p$ the normal, $\mathbf{n}(p)$, is vertical.  As a consequence, for $\epsilon=\epsilon(p)$ sufficiently small, near $p$, $\Sigma$ is the graph of a function over the disk $D_{\epsilon}(p)\subset \set{x_3=t_0}$.  Equivalently, let $C_{\epsilon}(p)=\set{(q,t)\in \Real^3 : q\in D_{\epsilon}}$ be the vertical cylinder over $p$ and $\pi_p$ be the natural projection $\pi_p: C_\epsilon(p)\to D_{\epsilon}(p)$. Then there is a neighborhood $U(p)$ of $p$ in $\Sigma$ so that $\pi_p$ restricted to $U(p)$ is a diffeomorphism onto $D_{\epsilon}(p)$.
For $\epsilon$ small enough,  $U(p)$ behaves with respect to $x_3$ as in Lemma \ref{HarmonLevelLem}.   

Let $\Sigma^\pm(p), \Sigma^{\geq t}(p)\subset \Sigma$ denote the sets given by Lemma \ref{HarmonLevelLem} \eqref{HLLii}, \eqref{HLLiii}.
As $\pi_p (U(p)\cap \set{x_3=t_0})=D_{\epsilon}(p)\cap \sigma_{t_0} $, if we let $\Omega^\pm(p)=\Omega^\pm_{t_0}\cap D_{\epsilon}(p)$ then either $\pi_p(\Sigma^+(p)) =\Omega^+(p)$ or $\pi_p(\Sigma^-(p))=\Omega^+(p)$.  We claim the parity of the identification is determined by whether the normal points up or down at $p$.  
Indeed, if the normal to $\Sigma$ points up at $p$ then, for $t>t_0$, $\pi_p(\Sigma^{\geq t}(p))= \pi_p(\Omega_t^-)$. Letting $t\to t_0$ gives $\pi_p(\Sigma^+ (p))=  \Omega^-(p)$. Conversely, if the normal to $\Sigma$ points down at $p$ then, for $t>t_0$,  $\pi_p( \Sigma^{\geq t}(p))= \pi_p(\Omega_t^+)$.  Letting $t\to t_0$ gives $\pi_p(\Sigma^+(p))=\Omega^+(p)$.
We further claim that the identification at $p$ determines the identification at $q=I(p)$.  Indeed, the identification is reversed and so the normals point in opposite directions -- in particular $q\neq p$.  
Figure \ref{Levelsets} illustrates this for a simple critical point.

To verify this we suppose, without loss of generality, that $\Sigma^+(p)$ is identified with $\Omega^+(p)$ -- i.e. the normal at $p$ is down. 
We will show that $\Sigma^-(q)$ is then identified with $\Omega^+(q)$ -- i.e. the normal is up.  
To begin the argument, we find a planar domain $B(p)$, a connected component of either $\Omega^+_{t_0}$ or $\Omega^-_{t_0}$, so that
\begin{enumerate}
\item \label{BDefi}$p\in \partial B(p)$,
\item \label{BDefii}$\partial B(p)\subset \sigma_{t_0}$,
\item \label{BDefiii}$ D_{\epsilon}(p)\cap B(p)$ consists of exactly one connected component, $B_0(p)$.  
\end{enumerate}
We justify the existence of $B(p)$ as follows: there are at least four
curves emanating from $p$ in $\Sigma'\cap \set{x_3=t_0}$ while
$\partial \sigma_{t_0}$ consists of only two points. Hence, there is a
connected component, $A_1$, of $\set{x_3=t_0}\backslash \sigma_{t_0}$
with $p\in \partial A_1$ and $\partial A_1\subset \sigma_{t_0}$.
As $A_1$ satisfies \eqref{BDefi} and \eqref{BDefii}, if $A_1\cap
D_{\epsilon}(p)$ has one component we are done. If there is more than
one component, they all lie in either $\Omega^+(p)$ or $\Omega^-(p)$
and so cannot be adjacent.  Thus, there is a simple closed curve
$\tau_1$ through $p$ lying in the closure of $A_1$ which bounds a
(topological) disk $\tilde A_1\subset \set{x_3 = t_0}$ so that
$\tilde{A}_1$ meets both $\Omega^+(p)$ and $\Omega^-(p)$. In
particular, $\set{x_3=t}\backslash A_1$ has a connected component,
$A_2\subset \tilde{A}_1\backslash \sigma_{t_0}$,  so $p\in \partial
A_2$ and $\partial A_2 \subset \sigma_{t_0}$.  If $A_2$ is not the desired
set, the same method produces a set $A_3$ disjoint from $A_1\cup A_2$,
so $p\in\partial A_3$ and $\partial A_3\subset \sigma_{t_0}$. Proceeding
in this fashion, because there are only finitely many components of
$D_{\epsilon}(p)\backslash \sigma_{t_0}$, one must eventually find
$A_{i_0}=B(p)$ satisfying \eqref{BDefiii}.

 Without loss of generality, we suppose $B(p) \subset \Omega^-_{t_0}$.
 Label the components of $\partial B(p)$  as $\alpha^1(p), \ldots,
 \alpha^k(p)$ so $p\in \alpha^1(p)$, and let $\alpha^1_0(p)=D_\epsilon
 (p)\cap \partial B_0(p)$.  There exists a single connected component
 $\Sigma_0^B(p)\subset \Sigma^-(p)$ with $\partial
 \Sigma_0^B(p)=\alpha_0^1(p)$ and $\pi_p (\Sigma_0^B(p))=B_0(p)$.
Let $I(\alpha^i(p))=\alpha^i(q)$.  We  will show that the $\alpha^i(q)$ are the
boundary of a connected planar domain $B(q)\subset\Omega^-_{-t_0}$
that satisfies $\eqref{BDefi}, \eqref{BDefii}$ and $\eqref{BDefiii}$
(with $p$ and $t_0$ replaced by $q$ and $-t_0$).
Notice \emph{a priori} there need not be any planar domain
in $\Omega^-_{-t_0}$ with boundary the $\alpha^i(q)$.  We will be able
to construct such a domain using topological properties of $I$ and by solving an appropriate Plateau problem.
 %parameterized by almost injective maps $\tilde{\alpha}^i(p)$. 
\begin{figure}[h]
         \centerline{
           \scalebox{.5}{
             \input{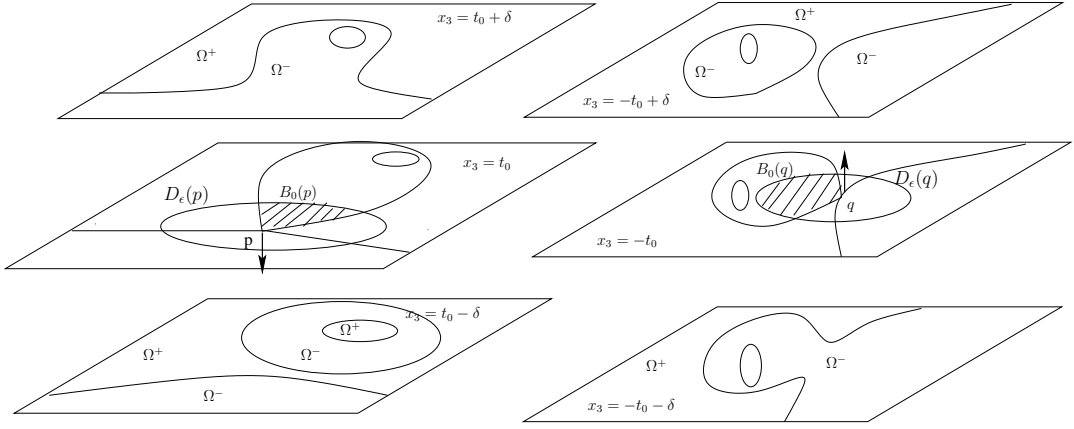}
             }
           }
         \caption{The left column shows level sets of $x_3$ near $p$.  The right shows the same near $q=I(p)$.
The shaded regions are $B_0(p)$ and $B_0(q)$.}
         \label{Levelsets}
       \end{figure}
Our argument exploits the existence of nice curves $\alpha^i_j(p)$
that are smoothly embedded is $\Sigma$ and pairwise disjoint (for
fixed $j$) and that converge to
$\alpha^i(p)$ in a $C^0$ sense and in the flat metric.
In particular, this allows us to think of the $\alpha^i(p)$ as cycles in $\Sigma$.
 The
existence of the $\alpha^i_j(p)$ follows from \eqref{HLLiv} of Lemma
\ref{HarmonLevelLem}.    As $I$ is a diffeomorphism we may then set $\alpha^i_j(q) =
I(\alpha^i_j(p))$ and obtain corresponding curves approximating
the $\alpha^i(q)$.  Using the $\alpha^i_j(p)$ and $ \alpha^i_j(q)$ together with the
maximum principle and the classification of surfaces we conclude that no
collection of either the $\alpha^i(p)$ or of the $\alpha^j(q)$ can be
null-homologous in $\Sigma$.  

% Away from critical points of $x_3$ the $\alpha^i(p)$ are smooth and
% embedded.  On the other hand, near each critical point $p_j
% \in \partial B(p)$,  the natural projection map $\pi_{p_j}$ and Lemma
% \ref{HarmonLevelLem} \eqref{HLLiv} give piecewise smooth
% parameterizations and allow us to think of each $\alpha^i(p)$ as a
% cycle in $\Sigma$.  Moreover, the $\alpha^i(p)$ are each the $C^0$
 %limit of a sequence of smooth embedded curves, which, in addition, converge in flat norm to $\alpha^i(p)$.
%Let $\alpha^i(q) :=I(\alpha^i(p))$ be piecewise smooth curves
%in $\{x_3=-t_0\}$ which we also may think of as cycles in $\Sigma$. 
%As $I$ is a diffeomorphism, the $\alpha^i(q)$ may
%also be approximated by smooth embedded curves. 

Recall the hyperelliptic
involution negates homology classes of $\Sigma$, thus 
\begin{equation}\label{HomologySum}
\sum_{i=1}^k
[\alpha^i(q)]=-\sum_{i=1}^k[\alpha^i(p)]
\end{equation}
 where $[\alpha^i(p)]$ and
$[\alpha^i(q)]$ denote the class in $H_1(\Sigma,\mathbb{Z})$ of the cycles
$\alpha^i(p)$ and $\alpha^i(q)$.  As
$\Sigma\backslash \Sigma'$ is an annulus, the inclusion map
induces an isomorphism between $H_1(\Sigma',\mathbb{Z})$ and
$H_1(\Sigma, \mathbb{Z})$. Inclusion also induces a map from
$H_1(\Sigma',\mathbb{Z})$ to $H_1( \bar{\Omega}^-,\mathbb{Z})$. It
follows that in $H_1(\bar{\Omega}^-,\mathbb{Z})$, $[\alpha^i(q)] $ and $[\alpha^i(p)]$ 
still satisfy \eqref{HomologySum}.  Thus, there is a
chain $C_0$ in $\bar{\Omega}^-$ with $\partial C_0=\sum_{i=1}^k\left(
  \alpha^i(p)+\alpha^i(q)\right)$.  Notice $\overline{B(p)}\subset
\bar{\Omega}^-$ may be thought of as a chain in $\bar{\Omega}^-$ with
$\partial \overline{B(p)}=\sum_{i=1}^k \alpha^i(p)$.  Hence,
$C_0-\overline{B(p)}$ is a chain in $\bar{\Omega}^-$ with $\partial (C_0-\overline{B(p)})=\sum_{i=1}^k \alpha^i(q)$.
Solving a constrained Plateau problem gives
a mass minimizing current $B'$ in $\Omega^-$ with $\partial B'=\sum_{i=1}^k
\alpha^i(q)$.  By the convex hull property, $B(q):=spt(B')\subset
\set{x_3=-t_0}$ and is a union of connected planar domains.  We claim
$B(q)$ is our desired domain.

For the sake of completeness we first discuss why
 $B'$ exists.   
As the $\alpha_j^i(q)$
converge in flat norm to $\alpha^i(q)$, for each $j$, the set of
curves are null-homologous in $\bar{\Omega}^-$.  Let
$B_j$ be the mass minimizing current in $\bar{\Omega}^-$ with
$\partial B_j=\sum_{i} \alpha^i_j(q)$; such a $B_j$ exists by direct
methods.  As $\bar{\Omega}^-$ is the complement of the union of open
sets with smooth boundary, Proposition 6.1 and Theorem 6.2 of
\cite{WhiteEllpEx} imply that $spt(B_j)\backslash
spt(\partial B_j)$ is a $C^{1, \alpha}$ surface in $\Real^3$ for some
$0<\alpha<1$.   Thus, $spt(B_j)$ is disjoint from the
singularities of $\partial \Omega^-$ and so we may apply the strong maximum principle of
Solomon and White \cite{SolWhite} to see that $spt(B_j)\backslash spt(\partial B_j)$ is either disjoint from
$\partial \Omega^-$ or a subset of $\Sigma$.  The latter case cannot
occur, for if it did the $\alpha_j^i(q)$, and hence the $\alpha^i(q)$, would be
null-homologous in $\Sigma$.  Thus, $spt(B_j) \backslash spt(\partial
B_j)\subset \Omega^-$.  We recover $B'$
from the $B_j$ by letting $j\to \infty$ and using standard compactness
theorems.

Let us now check that $B(q)$
is connected.  Denote by $\hat{B}(q)$ the connected
component of $B(q)$ with $\alpha^1(q)\subset \partial \hat{B}(q)$; if
$B(q)$ is not connected $\hat{B}(q)$ is a proper subset of $B(q)$.  In this
case, up to a relabeling,  $\partial \hat{B}(q)=\cup_{i=1}^{k'} \alpha^i(q)$ where
$k'<k$.  By the argument of the preceding two paragraphs, there is a
 a mass minimizing current
$B''$ in $\Omega^-$ with $\partial B''=\sum_{i=1}^{k'}\alpha^i(p)$.  As above, the convex hull property implies
$spt(B'')\subset \set{x_3=t_0}$.  This implies $B(p)$ is disconnected and so
is impossible.

By construction, $q\in \alpha^1(q)$ and $\partial B(q)=\cup_i \alpha^i(q)\subset \sigma_{-t_0} $.
Taking small enough values of $\epsilon$ (possibly differing at $p$ and $q$) we may ensure $U(q)\subset I(U(p))$.
Then $ \alpha^1_0(q)=\alpha^1(q)\cap U(q)\subset I(\alpha^1(p)\cap
U(p))=I(\alpha_0^1(p))$. As any component of $B(q)\cap
D_{\epsilon}(q)$ must have boundary containing $\alpha^1(q)\cap U(q)$
we conclude that  $B(q)\cap D_{\epsilon}(q)$ has only one component,
$B_0(q)$, satisfying  $ D_\epsilon(q)\cap \partial
B_0(q)=\alpha_0^1(q)$.  Hence, $B(q)$ is connected and satisfies
\eqref{BDefi}, \eqref{BDefii} and \eqref{BDefiii} as claimed.

Clearly, \[U(q)\cap \partial I(\Sigma_0^B(p)) =U(q)\cap I(\alpha_0^1(p)) =\alpha_0^1(q)\] and so \[\pi_q(  U(q)\cap I( \Sigma_0^B(p)) )=B_0(q)\subset \Omega^-(q).\]
However,  $I$ flips the sign of $x_3$ and so \[U(q) \cap
I(\Sigma_0^B(p))  \subset \Sigma^+(q).\]  Hence the identification at
$q$ is reverse what it is at $p$; that is the normal points up rather
than down. 
\end{proof}

\bibliographystyle{amsplain}
\bibliography{rotsym}
\end{document}